\documentclass[11pt,twoside,english]{amsart}
\advance\oddsidemargin by -1.0cm
\advance\evensidemargin by -1.0cm
\advance\topmargin by -1.0cm 
\usepackage{amssymb}
\usepackage{babel}
\usepackage{amstext}
\usepackage{amscd}   
\usepackage{epsfig}  
\usepackage{rotating}

\theoremstyle{plain}

\newtheorem{thm}{Theorem}[section]            
\newtheorem{prop}{Proposition}[section]
\theoremstyle{definition}

\newcommand{\bdm}{\begin{displaymath}}
\newcommand{\edm}{\end{displaymath}}
\newcommand{\be}{\begin{equation}}
\newcommand{\ee}{\end{equation}}
\newcommand{\ba}[1]{\begin{array}{#1}}
\newcommand{\ea}{\end{array}}

\newcommand{\btab}{\begin{tabular}}
\newcommand{\etab}{\end{tabular}}

\newcommand{\R}{\ensuremath{\mathbb{R}}}

\renewcommand{\S}{\ensuremath{\mathbb{S}}}

\newcommand{\G}{\ensuremath{\mathrm{G}}}
\newcommand{\J}{\ensuremath{\mathrm{J}}}

\newcommand{\SO}{\ensuremath{\mathrm{SO}}}

\begin{document}
\def\haken{\mathbin{\hbox to 6pt{%
                 \vrule height0.4pt width5pt depth0pt
                 \kern-.4pt
                 \vrule height6pt width0.4pt depth0pt\hss}}}
    \let \hook\intprod
\setcounter{equation}{0}
%
%
\thispagestyle{empty}
%
\date{\today}
\title[Nearly K\"ahler structures on $\S^6$]
{Nearly K\"ahler and nearly parallel $\G_2$-structures
on spheres}
%
%
%
\author{Thomas Friedrich}
\address{\hspace{-5mm} 
{\normalfont\ttfamily friedric@mathematik.hu-berlin.de}\newline
Institut f\"ur Mathematik \newline
Humboldt-Universit\"at zu Berlin\newline
Sitz: WBC Adlershof\newline
D-10099 Berlin, Germany}
%

\thanks{Supported by the SFB 647 "Raum, Zeit,Materie" and the SPP 1154 
``Globale Differentialgeometrie'' of the DFG}
\subjclass[2000]{Primary 53 C 25; Secondary 81 T 30}
\begin{abstract}
In some other context, the question was raised how many nearly K\"ahler
structures exist on the sphere $\S^6$ equipped with the 
standard Riemannian metric. In this short note, we prove that,
up to isometry, there exists only one.
 This is a consequence of the description
of the eigenspace to the eigenvalue $\lambda = 12$ 
of the Laplacian acting on $2$-forms. A similar result concerning
nearly parallel $\G_2$-structures on the round 
sphere $\S^7$ holds, too. An alternative proof by Riemannian Killing 
spinors is also indicated.
\end{abstract}
\maketitle
\pagestyle{headings}
%
%
%
\noindent
Consider the $6$-dimensional sphere $\S^6 \subset \R^7$ equipped with its
standard metric. Denote by $\Delta$ the Hodge-Laplace operator acting
an $2$-forms of $\S^6$ and consider the space
\bdm
E_{12} \ := \ \big\{ \omega^2 \in \Gamma(\Lambda^2(\S^6)) \, : \, d \, * \, 
\omega^2
\ = \ 0  , \ \ \Delta(\omega^2) \ = \ 12 \cdot \omega^2 \big\}  .
\edm 
This space is an $\SO(7)$-representation. Moreover, it 
coincides with the full eigenspace of the Laplace operator acting on
$2$-forms with eigenvalue $\lambda = 12$. 
\vspace{1mm}

\begin{prop}
The $\SO(7)$-representation $E_{12}$ is isomorphic to $\Lambda^3(\R^7)$. More 
precisely, for any $2$-form $\omega^2 \in E_{12}$, there exists a unique
algebraic $3$-form $A \in \Lambda^3(\R^7)$ such that
\bdm
\omega^2_x( y , \, z) \ = \ A(x , \, y , \, z)
\edm
holds at any point $x \in \S^6$  for any two tangent vectors
$y,z \in T_x(\S^6)$.
\end{prop}
\begin{proof}
It is easy to check that any $2$-form $\omega^2$ on $\S^6$
defined by  a $3$-form $A \in \Lambda^3(\R^7)$ as indicated
satisfies the differential equations $d * \omega^2 = 0 , \, 
\Delta(\omega^2) = 12 \cdot \omega^2$. Consequently, we obtain an
$\SO(7)$-equivariant map
\bdm
\Lambda^3(\R^7) \ \longrightarrow \, E_{12} .
\edm
Since $\Lambda^3(\R^7)$ is an irreducible $\SO(7)$-representation, the
map is injective. On the other hand, by Frobenius reciprocity, 
one computes the dimension of the eigenspace of the Laplace operator
on $2$-forms to the eigenvalue $\lambda = 12$. Its dimension equals $35$.
\end{proof}
\vspace{3mm}

\noindent
We recall some basic properties of nearly K\"ahler manifolds in dimension
six (see the paper \cite{AFS}). 
Let $(M^6, \, \J \, , \, g)$ be a nearly K\"ahler $6$-manifold.
Then it is an Einstein space with positive scalar curvature $\mathrm{Scal} >
0$. The K\"ahler form $\Omega$ satisfies the differential equations
\bdm
d \, * \, \Omega \ = \ 0 \, , \quad \Delta(\Omega) \ = \ \frac{2}{5} \cdot 
\mathrm{Scal} \cdot \Omega \, .
\edm
In particular, the K\"ahler form $\Omega^{\J}$ of {\it any} nearly K\"ahler
structure $(\S^6, \, \J  , \, g_{can})$ on the standard sphere 
$\S^6$ is a $2$-form on $\S^6$ satisfying
the equations $d \, * \, \Omega^{\J} =  0$ and $\Delta(\Omega^{\J}) = 
12 \cdot \Omega^{\J}$. 
This observation yields the following result.
\vspace{3mm}

\begin{prop}
The K\"ahler form $\Omega^{\J}$ of {\it any} nearly K\"ahler structure
 $(\S^6, \, \J , \, g_{can})$ on the standard sphere is given by an
algebraic $3$-form $A \in \Lambda^3(\R^7)$ via the formula
\bdm
\Omega_x^{\J}( y \, , \, z) \ = \ A(x \, , \, y \, , \, z)
\edm
where $x \in \S^6$ is a point in the sphere and $y,z \in T_x(\S^6)$ are
tangent vectors.
\end{prop}
\noindent
Since the K\"ahler form $\Omega^{\J}$ is a non-degenerate $2$-form  
at any point of
the sphere $\S^6$, the $3$-form $A \in \Lambda^3(\R^7)$ is a non-degenerate
vector cross product in the sense of Gray (see \cite{Brown}, \cite{Gray1}, 
\cite{Gray2}). For purely algebraic
reasons it follows that two forms of that type are equivalent
under the action of the group $\SO(7)$. Finally, we obtain the following
\vspace{3mm}

\begin{thm}
Let $(\S^6, \, \J , \, g_{can})$ be a nearly K\"ahler structure
on the standard $6$-sphere. Then the almost complex structure
$\J$ is conjugated -- under the action of the isometry group $\SO(7)$ --
to the standard nearly K\"ahler structure of $ \S^6$.
\end{thm}
\noindent
A similar argument applies in dimension seven, too.
\vspace{3mm}

\begin{thm}
Let $(\S^7, \, \omega , \, g_{can})$ be a nearly parallel $\G_2$-structure
on the standard $7$-sphere. Then it is
conjugated -- under the action of the isometry group $\SO(8)$ --
to the standard nearly parallel $\G_2$-structure of $ \S^7$.
\end{thm}
\vspace{3mm}

\noindent
{\bf Remark.}
Nearly K\"ahler structures in dimension six and nearly parallel structures
in dimension seven correspond to Riemannian Killing spinors. It is well-known
that the isometry group of the spheres $\S^6$ and $\S^7$ acts transitively
on the set of Killing spinor of length one. This observation yields a second 
proof of the latter Theorems (see \cite{FKMS} and \cite{Grun}).
\end{document}